\documentclass[11pt]{amsart}
\usepackage{amsmath}
\usepackage{amssymb}
\usepackage{amsthm}
\usepackage{color}

\usepackage{tikz}
\usetikzlibrary{calc,fadings,decorations.pathreplacing}

\newcommand\pgfmathsinandcos[3]{%
  \pgfmathsetmacro#1{sin(#3)}%
  \pgfmathsetmacro#2{cos(#3)}%
}
\newcommand\LongitudePlane[3][current plane]{%
  \pgfmathsinandcos\sinEl\cosEl{#2} 
  \pgfmathsinandcos\sint\cost{#3} 
  \tikzset{#1/.estyle={cm={\cost,\sint*\sinEl,0,\cosEl,(0,0)}}}
}
\newcommand\LatitudePlane[3][current plane]{%
  \pgfmathsinandcos\sinEl\cosEl{#2} 
  \pgfmathsinandcos\sint\cost{#3} 
  \pgfmathsetmacro\yshift{\cosEl*\sint}
  \tikzset{#1/.estyle={cm={\cost,0,0,\cost*\sinEl,(0,\yshift)}}} 
}
\newcommand\DrawLongitudeCircle[2][1]{
  \LongitudePlane{\angEl}{#2}
  \tikzset{current plane/.prefix style={scale=#1}}
  \pgfmathsetmacro\angVis{atan(sin(#2)*cos(\angEl)/sin(\angEl))} %
  \draw[current plane] (\angVis:1) arc (\angVis:\angVis+180:1);
  \draw[current plane,dashed] (\angVis-0:1) arc (\angVis-0:\angVis:1);
}
\newcommand\DrawLongitudeCircletwo[2][1]{
  \LongitudePlane{\angEl}{#2}
  \tikzset{current plane/.prefix style={scale=#1}}
  \pgfmathsetmacro\angVis{atan(sin(#2)*cos(\angEl)/sin(\angEl))} %
  \draw[current plane, dashed] (\angVis:1) arc (\angVis:\angVis+180:1);
  \draw[current plane,dashed] (\angVis-0:1) arc (\angVis-0:\angVis:1);
}

\newcommand\DrawLattwoCircle[2][1]{
  \LongitudePlane{\angEl}{#2}
  \tikzset{current plane/.prefix style={scale=#1}}
  \pgfmathsetmacro\angVis{atan(sin(#2)*cos(\angEl)/sin(\angEl))} %
  \draw[current plane, dashed] (\angVis:1) arc (\angVis:\angVis+180:1);
}

\newcommand\DrawLongCircle[2][1]{
  \LongitudePlane{\angEl}{#2}
  \tikzset{current plane/.prefix style={scale=#1}}
  \pgfmathsetmacro\angVis{atan(sin(#2)*cos(\angEl)/sin(\angEl))} %
  \draw[current plane,thin, dashed] (\angVis:1) arc (\angVis:\angVis+96:1);
}

\newcommand\DrawLongtwoCircle[2][1]{
  \LongitudePlane{\angEl}{#2}
  \tikzset{current plane/.prefix style={scale=#1}}
  \pgfmathsetmacro\angVis{atan(sin(#2)*cos(\angEl)/sin(\angEl))} %
  \draw[current plane,thin, dashed] (\angVis:1) arc (\angVis:\angVis+95:1);
}

\newcommand\DrawLatitudeCircle[2][1]{
  \LatitudePlane{\angEl}{#2}
  \tikzset{current plane/.prefix style={scale=#1}}
  \pgfmathsetmacro\sinVis{sin(#2)/cos(#2)*sin(\angEl)/cos(\angEl)}
  \pgfmathsetmacro\angVis{asin(min(1,max(\sinVis,-1)))}
  \draw[current plane, dashed] (\angVis:1) arc (\angVis:-\angVis-180:1);
}

\tikzset{%
  >=latex, 
  inner sep=0pt,%
  outer sep=2pt,%
  mark coordinate/.style={inner sep=0pt,outer sep=0pt,minimum size=3pt,
    fill=black,circle}%
}

\newtheorem{theorem}{Theorem}[section]
\newtheorem{corollary}[theorem]{Corollary}

\newtheorem{lemma}[theorem]{Lemma}

\theoremstyle{definition}

\theoremstyle{remark}

\def \t{\tilde}

\def \proj{\operatorname{P}}
\def \lacshadows{\operatorname{\sf Lsh}}
\def \lac{\operatorname{\sf Lac}}
\def \plac{\operatorname{\sf Lac_p}}
\def \pmax{\operatorname{\sf Max_p}}

\newtheorem{TheoC}{\textbf{\emph{Theorem C}}}
\newtheorem{TheoA}{\textbf{\emph{Theorem A}}}
\newtheorem{CoroB}{\textbf{\emph{Corollary B}}}

\newcommand{\dem}{\noindent {\it Proof. }}

\def\R{\mathbb{R}}

\def \la{\lambda}
\renewcommand{\le}{\leqslant}
\renewcommand{\ge}{\geqslant}

\newcommand{\fin}{\hspace*{\fill} $\square$ \vskip0.2cm}

\begin{document}

\newcommand{\demA}{\noindent {\bf Proof of Theorem A. }}
\renewcommand{\thefootnote}{}

\newcommand{\demB}{\noindent {\bf Proof of Theorem~B. }}

\title[On directional maximal operators]{Directional maximal operators and lacunarity in higher dimensions}


\address{Instituto de Ciencias Matem\'aticas CSIC-UAM-UC3M-UCM,
28049 Madrid, Spain}
\email{javier.parcet@icmat.es, keith.rogers@icmat.es}
\thanks{Supported by ERC grant 256997
 and MINECO grants MTM2010-16518, MTM2013-41780-P and SEV-2011-0087.
} 
\author{Javier Parcet}
\author{Keith M. Rogers}


\maketitle



\begin{abstract} We introduce a notion of lacunarity in higher dimensions for which we can
bound the associated directional maximal operators
 in $L^p(\R^n)$, with $p>1$. In particular, we are able to treat
the classes previously considered by
Nagel--Stein--Wainger, Sj\"ogren--Sj\"olin and Carbery. Closely related to this, we find a characterisation of the sets of directions which give rise to bounded maximal operators.
  The bounds enable Lebesgue-type differentiation of integrals in $L_{\text{loc}}^p(\mathbb{R}^n)$, replacing balls by tubes which point in these directions.
\end{abstract}


\vspace{1.5em}

\null\vspace{1.65em}
\section*{Introduction}

For $n \ge 2$ and a  set of directions $\Omega$ in the  unit sphere
$\mathbb{S}^{n-1}$, the  directional maximal operator $M_\Omega$ is
 defined, initially on Schwartz functions,~by
$$
M_\Omega f(x) =\sup_{\omega\in\Omega}\sup_{r>0} \frac{1}{2r}\int_{-r}^r |f(x-t\omega)|\,dt.
$$
If $\Omega$ consists of a single direction and $p>1$, the boundedness of $M_\Omega$ from $L^p(\R^n)$ to $L^p(\R^n)$ follows from the Hardy--Littlewood maximal theorem. 
If~$\Omega$ consists of many directions, two questions naturally arise:
\begin{itemize}
\!\!\item[(i)] When $\Omega$ is an arbitrary finite subset of $\mathbb{S}^{n-1}$, a fundamental problem is to determine the best bounds for  the $L^p$--operator norm of $M_{\Omega}$ as a function of $|\Omega|$ and $p$.

\vskip3pt

\item[(ii)] When $\Omega$ is an infinite subset of $\mathbb{S}^{n-1}$, one can also ask for  conditions on the directions which ensure that $M_{\Omega}$ is $L^p$--bounded.
\end{itemize}  

In two dimensions, the questions have been answered with remarkable
accuracy (see \cite{Co,Str,K0,K} for the first question, \cite{Str0,
CF, NSW, SS, B} for the second question, or \cite{ASV,ASV2, A, KL} which address
the two questions in a unified way), however much less is known in higher dimensions (see \cite{W,Bo2,KLT} for  the first
question and \cite{NSW, SS, C} for the second). 
We will prove a localization principle -- subsets of the directions can be considered independently from the rest of the directions --  from which we draw conclusions for the second question in three dimensions and more. A fundamental difference between the two--dimensional problem, with directions in $\mathbb{S}^1$, and that of higher dimensions is that we are no longer able to order the directions.

We partition the unit ball in such a way that it resembles a peeled orange with infinitely many segments. 
In three dimensions, we make three partitions, each time with a different axis of partition. The independent sets of our localisation principle will be contained in these segments.
More precisely and more generally, for $\sigma\in\Sigma$, where $$\Sigma \equiv\Sigma(n)= \{(j,k) :  1 \le j < k \le n\}$$ we
consider $\{\theta_{\sigma\!,i}\}_{i\in \mathbb{Z}}$ that
satisfy $0< \theta_{\sigma\!, i+1}\le \lambda_{\sigma}\,
\theta_{\sigma\!,i}$ with lacunary constants $0<\lambda_{\sigma}<1$.
Then, for an orthonormal basis $(e_1,\ldots,e_n)$, we divide the
directions into the subsets $\Omega_{\sigma\!,i}$ defined by 
$$ \Omega_{\sigma\!,i}=\Big\{\, \omega \in \Omega \ :\ \theta_{\sigma\!,
i+1}< \Big|\frac{\omega\cdot e_k}{\omega\cdot e_j}\Big|\le \theta_{\sigma\!, i}\, \Big\}$$
(see Figure 1). Note that the segments become thinner as $i$ converges to $\pm \infty$ and the partition of $\Omega$ is completed by including the set $\Omega_{\sigma\!,\infty}=\Omega\cap (e_j^\perp\cup e_k^\perp)$.

Writing $\mathbb{Z}^*=\mathbb{Z}\cup\{\infty\}$, we prove the
following localisation principle which recalls the separation of dyadic frequency scales provided by Littlewood--Paley theory (see also \cite{NT,CS} for another kind of one--dimensional localisation). A difference is that we have many lacunary partitions instead of one, however the result is sharp in the sense that the supremum over
partitions must be taken over the whole of $\Sigma$. Nor could it be made more flexible by allowing the segments to \lq accumulate' away from the hyperplanes perpendicular to the basis vectors.

\begin{TheoA}\label{ours2} \emph{Let $n\ge 2$ and $p>1$. Then
\begin{align*}
 \|M_{\Omega}\|_{p\to p} \ \le\  C\, \sup_{\sigma \in \Sigma}\, \sup_{i\in\mathbb{Z}^*}\|M_{\Omega_{\sigma\!,i}}\|_{p\to p},
\end{align*}
where $C$ depends only on $n$, $p$ and the lacunary constants $\lambda_\sigma$ for $\sigma \in \Sigma$.}
\end{TheoA}

As with the almost orthogonality principle of Alfonseca, Soria and Vargas  in
two dimensions~\cite{ASV, ASV2, A}, we recover the previously known results for 
question~(ii) in higher dimensions.
 Nagel, Stein and
Wainger \cite{NSW} proved the $L^p$--boundedness of the maximal operator associated to the
directions $$\{(\vartheta_i^{a_1},\ldots,\vartheta_i^{a_n})\}_{i\ge 1},$$
where $0<a_1<\ldots<a_n$ and $0< \vartheta_{i+1}\le \lambda\, \vartheta_{i}$ with lacunary constant
$0<\lambda<1$. We can apply Theorem~A
with~$\theta_{\sigma\!,i}=
\vartheta_i^{a_k-a_j}$ and $\la_{\sigma}=\lambda^{a_k-a_j}$, where $\sigma=(j,k)$, reducing the problem to that of a
single direction. Note that it makes no difference if the directions are normalised to live on the unit sphere or not.
 On the other hand, Carbery~\cite{C} proved that the maximal operator associated to the directions \begin{equation}\label{carbery}\{(2^{k_1},\ldots,2^{k_n})\}_{k_1,\ldots,k_n\in\mathbb{Z}}\end{equation} is $L^p$--bounded with $p>1$. Taking $\theta_{\sigma\!,i} = 2^{-i}$, the resulting sets of directions $\Omega_{\sigma\!,i}$ are restricted to  $(n-1)$--dimensional hyperplanes, so that by choosing a suitable basis and applying Fubini's theorem, we reduce to the $(n-1)$--dimensional problem. Iterating the process, applying the $d$--dimensional version of Theorem~A if the directions are restricted to a $d$--dimensional hyperplane, we eventually end up with isolated directions as before.

\begin{figure}
\centering
\begin{tikzpicture} 
\def\Ra{3} 
\def\angEl{25} 
\filldraw[ball color=white] (0,0) circle (\Ra);
\foreach \t in {45,-45,-43.21008939, -41.42366563, -37.87498366, -30.96375654, -18.43494883, 0.,18.43494883,30.96375655,37.87498363,41.42366563,43.21008936,46.78991061, 48.57633437, 52.12501634, 59.03624346, 71.56505117, 90., 108.4349488, 120.9637566, 127.8749836, 131.4236656, 133.2100894} { \DrawLongitudeCircle[\Ra]{\t} }
\DrawLatitudeCircle[\Ra]{0}

\pgfmathsetmacro\H{\Ra*cos(\angEl)}
\coordinate[mark coordinate] (N) at (0,\H);
\node[above=-190pt] at (N) {{\bf Figure 1.} A partition for $\sigma=(1,2)$ of a three--dimensional dissection.};
\pgfmathsetmacro\H{\Ra*cos(\angEl)}
\coordinate[mark coordinate] (N) at (0,\H);
\node[above=11pt] at (N) {$e_{3}$};
\node[above=1pt] at (N) {$\uparrow$};
\coordinate[mark coordinate] (eone) at (-2.13,\H-3.61);
\node[above=-0.5pt, left=-4pt] at (eone) {$\leftarrow$};
\node[left=6pt] at (eone) {$e_1$};
\coordinate[mark coordinate] (etwo) at (2.13,\H-3.61);
\node[above=-0.5pt, right=-4pt] at (etwo) {$\rightarrow$};
\node[right=7pt] at (etwo) {$e_2$};

\end{tikzpicture}
\end{figure}

It is not sufficient to constrain the angles
between an infinite number of directions if they are to give rise to a bounded maximal
operator in
higher dimensions, even if the directions live inside certain (indeed most) smooth curves. However Theorem A suggests a definition of lacunarity that
gives rise to bounded maximal operators in general. Given $\Omega\subset\mathbb{S}^{n-1}$, an orthonormal
basis of $\text{span}(\Omega)=\R^d$ with $d\le n$, and lacunary sequences
$\{\theta_{\sigma\!,i}\}_{i\in\mathbb{Z}}$, define
 partitions
$\{\Omega_{\sigma\!,i}\}_{i\in\mathbb{Z}^*}$ for each $\sigma\in \Sigma(d)$. We call such a choice of $\frac{1}{2}d(d-1)$
partitions a {\it dissection}. We say that $\Omega$ is
\begin{itemize}
\item  {\it  lacunary of order~0} if it consists of a single direction
\item  {\it lacunary of order $
L$}  if there is a dissection for which the sets~$\Omega_{\sigma\!,i}$ are lacunary of order $\le L-1$ for all $i\in
\mathbb{Z}^*$ and $\sigma \in \Sigma(d)$, with
uniformly bounded lacunary constants.
\end{itemize}
We say that a set of directions is {\it lacunary} if it is a finite union of sets which are lacunary of finite order and we denote the class of such sets by $\lac(n)$. Note that the class contains sets of directions which are confined to $d$--dimensional subspaces with $1\le d\le n$.

 According to this definition, the Nagel--Stein--Wainger directions
are lacunary of order $1$ and the Carbery directions are lacunary
of order $n-1$. 
By repeatedly applying Theorem~A as before, if $\Omega$ is lacunary (of finite order), then  $M_\Omega$ is
$L^p(\R^n)$--bounded with $p>1$. This extends the two--dimensional
result due to Sj\"ogren--Sj\"olin~\cite{SS} (the union of~$K$ sets
of directions of lacunary order $L$ with respect to their
definition, is lacunary of order $
2KL+1$ with respect to ours). We have broken with the two--dimensional tradition (which is one--dimensional in the sense that the directions are contained in a circle), whereby the label \lq lacunary' is reserved for the sets which are lacunary of order one. This is because such sets are less special in higher dimensions (they cannot necessarily be represented as a sequence for example) and at first glance the set of directions defined in \eqref{carbery} seems as deserving of the label \lq lacunary' as any other.

As a corollary we obtain a generalisation of the Fundamental Theorem of Calculus. After a suitably fine finite splitting of the directions, the operator
$M_\Omega$ can be composed with one--dimensional Hardy--Littlewood
maximal operators to dominate a constant multiple of the maximal
operator $\mathcal{M}_\Omega$ defined by
$$
\mathcal{M}_\Omega f(x) = \sup_{x \in T \in \mathcal{T}_\Omega}
\frac{1}{|T|}\int_{T} |f(y)|\,dy.
$$
Here,
$\mathcal{T}_\Omega$ denotes the family of tubes which point in a
direction of~$\Omega$. Standard density arguments yield the following Lebesgue
type differentiation result. 

\begin{CoroB} \emph{Let $n\ge 2$ and  $\Omega\in \lac(n)$. Then
 \begin{equation*}\label{conv}
\lim_{\substack{x\in T\in \mathcal{T}_{\Omega}\\
\mathrm{diam}(T) \to 0}} \frac{1}{|T|}\!\int_{T} f(y)\,dy=f(x),
\quad \mbox{a.e.} \ x\in \R^n,
\end{equation*}
for all $f \in L_{\mathrm{loc}}^p(\mathbb{R}^n)$ with $p>1$.}   
\end{CoroB}

Sets of directions which give rise to unbounded $M_\Omega$ and $\mathcal{M}_\Omega$ can be considered if we place further restrictions on the tubes $\mathcal{T}_\Omega$. Most commonly the eccentricity (length/width) of the tubes is fixed, especially when treating question (i) above. For question (ii), C\'ordoba~\cite{C3} proved that the associated maximal operator is bounded, with a logarithmic dependency on the eccentricity,  if the directions are restricted to a curve which intersects the hyperplanes of $\R^n$ no more than  a uniformly bounded number of times. 

We turn now to the question of characterising the sets of directions $\Omega$ for which $M_\Omega$ is bounded from $L^p(\R^n)$ to $L^p(\R^n)$. We denote the class of such sets by $\pmax(n)$. Bateman \cite{B} proved that $\pmax (2)\subset \lac(2)$ which, combined with the result of Sj\"ogren--Sj\"olin, yields the equivalence
$$\pmax (2)\equiv \lac(2),\quad 1<p<\infty.$$
We do not know if
this is true in higher dimensions, however we characterise $\pmax(n)$ using a
formally larger class. For this we take advantage of a quantitive version of Bateman's theorem via projections onto two--dimensional subspaces $\Pi\subset \R^n$. Given  a set of directions $\Omega$, we
define the {\it shadow of $\Omega$ on~$\,\Pi\,$} to be the normalised orthogonal projection onto $\Pi$, so that the shadow lives in a copy of $\mathbb{S}^1\subset \Pi$.
We denote by $\lacshadows(n)$ the class of sets of directions whose
shadows are all lacunary, where the lacunary constants are bounded uniformly away from one, and the number of sets in the finite unions and the lacunary orders are uniformly bounded above.

In order to define our characterising class, we fix an auxiliary $\varepsilon>0$ and suppose that $\text{span}(\Omega)=\R^d$.  We say that a set
$\Omega_{\widetilde{\sigma}\!,i_{\widetilde{\sigma}}}$ of a dissection is {\it dominating} if it satisfies
$$
\|M_{\Omega_{\sigma\!,i}}\|_{p\to p}\le
\|M_{\Omega_{\widetilde{\sigma}\!,i_{\widetilde{\sigma}}}}\|_{p\to p}+\varepsilon\quad \text{for all}\quad
i\in \mathbb{Z}^*, \ \sigma \in \Sigma(d).
$$ 
Similarly to before, we say that $\Omega$ is
\begin{itemize}
\item  {\it  $p$--lacunary of order~0} if it consists of a single direction
\item {\it $p$--lacunary of order $
L$}  if there is a dissection with
a dominating set  which is
$p$--lacunary of order $\le L-1$. 
\end{itemize}
We say that a set of directions is {\it $p$--lacunary} if it is a finite union of sets which are $p$--lacunary of finite order and we denote the class of such sets by $\plac(n,\varepsilon)$. Finally we write
$$
\plac(n):=\bigcap_{\varepsilon>0} \plac(n,\varepsilon).
$$

In the following
equivalence we see that the directions which give rise to bounded
maximal operators can be no worse, loosely speaking, than directions
that can be divided into isolated directions by a finite number of
lacunary dissections.

\begin{TheoC} \emph{Let $n\ge 2$ and $1<p<\infty$. Then
$$
\lac(n)\, \subset\, \pmax(n)\equiv \plac(n)\, \subset\, \lacshadows(n).
$$}
\end{TheoC}

 With $n=2$, these classes coincide of course, and so $\pmax(2)$ is the same for all $1<p<\infty$. It is tempting to
suppose that this is also true in higher dimensions -- it seems reasonable to expect that a member of $\lacshadows(n)$ could be dissected into isolated directions -- however it may
also be  that $\pmax(n)$ grows with $p$.  In any case, given the nature of the definitions of $\lac(n)$ and $\lacshadows(n)$, we
see that $\pmax(n)$ is not so far from a purely
two--dimensional concept. This should be compared with \cite{KLT2} (see also \cite{G}),
where Kakeya sets in $\R^3$ with Minkowski dimension sufficiently close to 5/2 were shown
to have a \lq planiness' property.

In the following section,  we prove Theorem A which implies the first inclusion of  Theorem C. The key ingredient is a nonlinear and nonpositive partition of a hyperplane in which we compensate for the points which are covered more than once by removing smaller sets. That is to say, our partition of unity is more like a covering (normally completely inadequate on the frequency side), but by adding and subtracting enough  times we are able to partition the hyperplane with intersections of tensor products of two--dimensional cones. These give rise to {\it a priori} frightening nonlinear terms, however they are dealt with later in a reasonably trivial fashion. In the
second section we prove the equivalence and the final inclusion of Theorem~C. Unusually in this context, this follows by a topological argument. In the third section, we justify a number of remarks from
above by constructing sets of apparently well--behaved directions for which the associated maximal operators are unbounded. In the final section, we provide a corollary for the maximal directional Hilbert transform. Some of these results  were announced in \cite{PR}.

\section{Proof of Theorem A}

By a finite splitting we can suppose that the
directions $\Omega$ are contained in the first open \lq octant' of the
unit sphere $\mathbb{S}^{n-1}\cap \mathbb{R}_{+}^n$. We consider
intersections of the segments  to obtain cells of directions
$$\Omega_{\bf i} = \bigcap_{\sigma \in \Sigma} \Omega_{\sigma\!,
i_{\!\sigma}}\quad \text{for each}\quad {\bf i}=(i_\sigma)_{\sigma \in \Sigma}
\in\mathbb{Z}^{\Sigma}.$$  This yields a finer partition than those
of the introduction;
$$\Omega \, = \, \bigcup_{{\bf i} \in \mathbb{Z}^{\Sigma}} \Omega_{\bf i}\qquad \text{so that}\qquad M_\Omega=\sup_{{\bf i} \in \mathbb{Z}^{\Sigma}} M_{\Omega_{\bf i}}.$$ Note that many of the cells are empty, however
we will see that this overdetermination is somehow unavoidable. Let
$K_{\sigma\!,i}$ denote the convolution operator associated to a
Fourier multiplier $\psi_{\sigma\!,i}$, smooth on $\mathbb{R}^n\backslash\{0\}$, equal to one on
$$
\Psi_{\!\sigma\!,i}=\Big\{\, \xi\in \mathbb{R}^n\ :\ \frac{1}{n}\theta_{\sigma\!,
i+1}< -\frac{\xi_j}{\xi_k}\le n\theta_{\sigma\!, i}\,\Big\},$$
and supported in a similar cone with $n$ replaced by $n+1$.  

The key geometric fact used in the proof of the following lemma is
that the hyperplane perpendicular
 to $\omega$ is contained in $\cup_{\sigma\in\Sigma} \Psi_{\!\sigma\!,i_{\!\sigma}}$ for all $\omega\in \Omega_{\bf i}$.  This is no longer true if, in the definition of the cones,  $n$ is replaced by a constant strictly less than $n-1$.  At this point we do not use that the dividing sequences are lacunary.

\begin{lemma} \label{key0} Let $p>1$. Then 
\begin{equation*}\label{corp}
\|M_\Omega\|_{p\to p} \, \le \ C\!\sup_{\emptyset\neq \Gamma \subset
\Sigma}\Big\| \sup_{{\bf i} \in \mathbb{Z}^{\Sigma}} M_{\Omega_{\bf
i}} \prod_{\sigma \in \Gamma} K_{\sigma\!, i_{\!\sigma}}\Big\|_{p\to
p},
\end{equation*}
where $C$ depends only on $n$ and $p$.
\end{lemma}

\dem
Fix a nonnegative, even, smooth function $m_{\text{o}}^\vee$ which
is positive on $[-1,1]$ and with sufficient decay
so that, for positive functions, $M_\Omega
f$ is pointwise equivalent to
$$
\sup_{\omega\in\Omega}\sup_{r>0}\Big|\frac{1}{r} \int
m_{\text{o}}^\vee(\tfrac{t}{r})f(\,\cdot\,-t\omega)\,dt\Big|.
$$
As the operator norm of $M_\Omega$ can be realised by testing on positive functions, we can work with the maximal operator $\widetilde{M}_\Omega$ defined by
$$
f\mapsto \sup_{\omega\in\Omega}\sup_{r>0}\Big|\frac{1}{r} \int
m_{\text{o}}^\vee(\tfrac{t}{r})f(\,\cdot\,-t\omega)\,dt\Big|,
$$
which is more amenable to Fourier analysis.
Throughout,  $^\wedge$ and $^\vee$ denote the
Fourier transform and inverse transform, respectively. One can
calculate that
$$\Big(\frac{1}{r} \int m_{\text{o}}^\vee(\tfrac{t}{r})f(\,\cdot\,-t\omega)\,dt\Big)^\wedge(\xi)=m_{\text{o}}(r \omega\cdot\xi)f^\wedge(\xi).
$$
It will simplify things to take $m_{\text{o}}$ supported in
$[-1,1]$, which can be arranged by choosing
$m_{\text{o}}=\phi_\text{o}\ast\phi_\text{o}$ where $\phi_\text{o}$
is an even, smooth function supported in $[-1/2,1/2]$. We also fix
a smooth function $\eta_\text{o}$, supported in the ball of radius
$4n^2$, centred at the origin and equal to one on the concentric ball
of radius $2n^2$, and consider the operator
$$
f\mapsto\sup_{\omega\in\Omega}\sup_{r>0}  |S_{r\!,\omega}f|,
$$
where
$
\big(S_{r\!,\omega}f\big)^\wedge(\xi)=\eta_\text{o}\big(r(\omega_1\xi_1,\ldots,\omega_n\xi_n)\big)m_{\text{o}}(r
\omega\cdot\xi)f^\wedge(\xi).
$
This is pointwise dominated  by a constant multiple of the
strong maximal operator~$\mathcal{M}_{\text{str}}$, which can be
bounded by iterated applications of the one--dimensional
Hardy--Littlewood maximal theorem. Defining
$m$ by $$m(\xi)=(1-\eta_\text{o})(\xi)m_{\text{o}}({\bf 1}\cdot\xi)$$
with ${\bf 1}=(1,\ldots,1)$, we are left with the maximal operator
$T_\Omega$ defined by
$$
 f\mapsto\sup_{\omega\in\Omega}\sup_{r>0}\big|T_{r\!,\omega}f|,
$$
where $ (T_{r\!,\omega} f)^\wedge(\xi) =
m\big(r(\omega_1\xi_1,\ldots,\omega_n\xi_n)\big)f^\wedge(\xi).$ A variant of this reduction was originally employed by Nagel, Stein and Wainger \cite{NSW}. 

It will suffice to prove the pointwise estimate
\begin{equation}\label{lk}
T_\Omega f \ \le \ \sum_{\emptyset\neq \Gamma \subset \Sigma}\,
\sup_{{\bf i} \in \mathbb{Z}^{\Sigma}}\, T_{\Omega_{\bf i}}
\Big[\prod_{\sigma \in \Gamma} K_{\sigma\!, i_{\!\sigma}}\Big]f.
\end{equation}
The desired $L^p$--estimate then follows by combining with the inequalities
$$
\widetilde{M}_\Omega f\le C(\mathcal{M}_{\text{str}}f +T_{\Omega} f),\quad T_{\Omega_{{\bf i}}} f\le
C\big(\mathcal{M}_{\text{str}}f+M_{\Omega_{{\bf i}}}f\big),
$$
and, when $\Omega_{{\bf i}}\neq \emptyset$,
\begin{align*}
\Big\| \sup_{{\bf i} \in \mathbb{Z}^{\Sigma}}
\mathcal{M}_{\text{str}} \prod_{\sigma \in \Gamma} K_{\sigma\!,
i_{\!\sigma}}\Big\|_{p\to p} &\le C\,\Big\| \sup_{{\bf i} \in
\mathbb{Z}^{\Sigma}} M_{\Omega_{{\bf i}}} \prod_{\sigma \in \Gamma}
K_{\sigma\!, i_{\!\sigma}}\Big\|_{p\to p}.
\end{align*}
The final inequality is a trivial consequence of the boundedness of the strong maximal operator, combined with the fact that $|f|\le M_{\Omega_{{\bf i}}} f$.

Before proving \eqref{lk}, we motivate why it is reasonable to hope that it should be true. As suggested earlier,  the frequency support of $T_{r\!,\omega}f$ is contained in the union of $\Psi_{\!\sigma\!,i_{\!\sigma}}$ whenever $\omega\in \Omega_{\bf i}$ with ${\bf i} = (i_\sigma)_{\sigma \in \Sigma}$. If this covering were in fact a partition, we would obtain
$$
T_{r\!,\omega}f \ = \
\sum_{\sigma\in \Sigma}^{\null} 
\, T_{r\!,\omega} 
K_{\sigma\!,i_{\!\sigma}}f,\qquad \omega\in \Omega_{\bf i},
$$
and so, recalling that $T_\Omega f = \sup_{{\bf i} \in \mathbb{Z}^{\Sigma}} T_{\Omega_{{\bf i}}}f$, a simplified version of \eqref{lk}, with less terms on the right-hand side, would follow easily.  Now the conic supports do not form a partition and so to compensate we remove the pairwise intersections of the cones and then add back the intersections of each triple of cones, and so on, until we  obtain a partition. Some of these intersections may in fact be empty, but we ignore this as there is no advantage for us to have less terms in the sum. Indeed, we will see that for our purposes there is no difference between the earlier simplified version and the following complicated looking formula. In three dimensions we can identify $\sigma=(1,2),(1,3),(2,3)$ with $3,2,1$, respectively, and the process yields the identity
$$
T_{r\!,\omega}f \ = \ \sum_{1\le j\le 3}^{\null} 
 T_{r\!,\omega} 
K_{j,i_{j}}f\ - \!\!\!\sum_{1\le j<k\le3}^{\null} 
\, T_{r\!,\omega} 
K_{j,i_{j}}K_{k,i_{k}}f\, + T_{r\!,\omega} 
K_{1,i_{1}}K_{2,i_{2}}K_{3,i_{3}}f
$$
plus a remainder term.
More generally, we obtain
$$
T_{r\!,\omega}f \ = \ 
\sum_{\emptyset\neq \Gamma \subset \Sigma}^{\null} (-1)^{|\Gamma|+1}
\, T_{r\!,\omega} \Big[\prod_{\sigma \in \Gamma}
K_{\sigma\!,i_{\!\sigma}}\Big]f\,+\,T_{r\!,\omega}\mathcal{R}_{{\bf i}}f. 
$$
 In effect, we have expanded the polynomial $1 - \prod_\sigma (1 - x_\sigma)$, and so 
the remainder $\mathcal{R}_{{\bf i}}$ is given by $$(\mathcal{R}_{{\bf
i}}f)^\wedge(\xi) = \prod_{\sigma \in
\Sigma}(1-\psi_{\sigma\!,i_{\!\sigma}})(\xi) f^\wedge(\xi).$$

In contrast with the operators $K_{\sigma\!,i}$, which are essentially
two--dimensional, the operators $\mathcal{R}_{\bf i}$ are genuinely
higher--dimensional objects, however once we see that the
multiplier associated to $T_{r\!,\omega}\mathcal{R}_{{\bf i}}$ is
identically zero whenever $\omega \in \Omega_{\bf i}$ and $r>0$,
\begin{equation}\label{empt}
m\big(r(\omega_1\xi_1,\ldots,\omega_n\xi_n)\big)\prod_{\sigma \in
\Sigma}(1-\psi_{\sigma\!,i_{\!\sigma}})(\xi) \, \equiv \, 0,
\end{equation}
we obtain 
$$
T_{r\!,\omega}f \ = \ 
\sum_{\emptyset\neq \Gamma \subset \Sigma}^{\null} (-1)^{|\Gamma|+1}
\, T_{r\!,\omega} \Big[\prod_{\sigma \in \Gamma}
K_{\sigma\!,i_{\!\sigma}}\Big]f\, ,\qquad \omega\in \Omega_{\bf i} 
$$
which yields 
\eqref{lk}. Given that the cones are invariant under scaling, by taking $r$ large, \eqref{empt} is little more than the assertion that the hyperplane is covered by the cones.

After the scaling
$\omega_j\xi_j\to \xi_j$ for $1 \le j \le n$, it will suffice to
prove that the region defined by
\begin{equation} \label{Cond1}
\Big| \sum_{j=1}^n \xi_j \Big| \le \frac{1}{r} \qquad \mbox{and}
\qquad \Big( \sum_{j=1}^n \xi_j^2 \Big)^{1/2} \ge \frac{2n^2}{r}
\end{equation}
and
\begin{equation*} \label{Cond2}
 -\frac{\xi_j}{\xi_k}\frac{\omega_k}{\omega_j}\le \frac{1}{n}\theta_{\sigma\!,
i+1}  \ \ \textrm{or}\ \ -\frac{\xi_j}{\xi_k}\frac{\omega_k}{\omega_j}>  n\theta_{\sigma\!, i}  \quad \mbox{for all} \quad
\sigma \in \Sigma
\end{equation*}
is empty. As $\omega\in \Omega_{\bf i}$, we see that the complements of the scaled cones are contained in
\begin{equation}\label{Cond3}
 -\frac{\xi_j}{\xi_k}< \frac{1}{n}  \quad \textrm{or}\quad -\frac{\xi_j}{\xi_k}> n  \quad \mbox{for all} \quad
\sigma \in \Sigma.
\end{equation}
We suppose for a contradiction that the region defined by
\eqref{Cond1} and \eqref{Cond3} is not empty. It is clear by
comparing the inequalities in \eqref{Cond1} that the components of a
vector $\xi$ in this region cannot all have the same sign. By
symmetric invariance of the conditions, we may suppose that
$$\xi_1,\ldots,\xi_{m-1}\ge 0\qquad \text{and} \qquad \xi_m,\ldots,\xi_{n}< 0$$ for
some $1 < m \le n$. We can also suppose without loss of generality
that $|\xi_1|\ge|\xi_j|$ for all $j>1$ and $|\xi_m|\ge|\xi_j|$ for
all $j>m$. Then taking $j=1$ and $k=m$ in~\eqref{Cond3} we see that $|\xi_1|\ge n|\xi_m|$. On the other hand, by the first
condition of \eqref{Cond1},
\begin{align*}
|\xi_1|-(n-1)|\xi_m| \le \Big| \sum_{j=1}^n \xi_j \Big| \le \frac{1}{r}.
\end{align*}
Combining the two estimates we obtain $ |\xi_1|\le n/r.$ Since
$|\xi_1| \ge |\xi_j|$ for $j > 1$, this yields $$
|\xi_1|+\ldots+|\xi_n|\le \frac{n^2}{r} $$ which contradicts the second
inequality in \eqref{Cond1}. Thus, $T_{r\!,\omega}\mathcal{R}_{\bf i}
\equiv 0$  whenever $r>0$ and $\omega \in \Omega_{\bf i}$, and we
are done. \fin

\vspace{0.5em}

\noindent We will also require the following square function estimates which follow easily from the two--dimensional theory.

\begin{lemma}\label{square} Let $1<p<\infty$ and $\Gamma \subset \Sigma$. Then
\begin{equation*}
\Big\|\Big(\sum_{{\bf i} \in \mathbb{Z}^\Gamma} \big|
\Big[\prod_{\sigma \in \Gamma} K_{\sigma\!,i_{\!\sigma}}\Big] f
\big|^2 \Big)^{\frac12} \Big\|_{{p}}  \le\, C\,\|f\|_p,\quad
\text{where}\quad {\bf i} = (i_\sigma)_{\sigma \in \Gamma},
\end{equation*}
and $C$ depends only on $|\Gamma|$, $p$ and the lacunary constants
$\lambda_\sigma$.
\end{lemma}

\dem In order to bound directional maximal operators in  $L^2$,  the required square
function estimate, with $p=2$,  follows directly from Plancherel's theorem and the
finite overlapping of the supports of $\{\psi_{\sigma\!,i}\}_{i\in
\mathbb{Z}}$. This is where we use the lacunarity of the sequences
~$\{\theta_{\sigma\!,i}\}_{i\in \mathbb{Z}}$. When $p\neq 2$, by a
standard randomisation argument, using Khintchine's inequality, the
square function estimates follow from the uniform
$L^p$--boundedness, independent of the choice of the signs, of the
Fourier multiplier operators
$$
f\mapsto  \Big(\sum_{{\bf i} \in \mathbb{Z}^\Gamma} \pm
\prod_{\sigma\in \Gamma}
\psi_{\sigma\!,i_{\!\sigma}}\,f^\wedge\Big)^\vee.
$$
This in turn is a consequence of the Marcinkiewicz multiplier
theorem (see for example \cite[pp. 109]{Stein}), for which it suffices to check a number of conditions
involving integrals of derivatives of the multipliers. After
applying the product rule, the calculation reduces to the
case~$|\Gamma|=1$. Applying Fubini's theorem so as to ignore the
trivial variables, this was originally checked by A. C\'ordoba and
R. Fefferman~\cite[Section 4]{CF0}  in their proof of a
two--dimensional angular Littlewood--Paley inequality. Again,  the calculation relies on the  lacunarity of the dividing sequences. \fin

\vspace{0.5em}

Armed with these lemmas, the proof is completed easily as follows. In order to establish the idea, we treat the easiest case first.

\vspace{1em}

\noindent\textbf{Case $n=3$ and $p=2$.} We can identify $\sigma=(1,2),(1,3),(2,3)$ with $3,2,1$, respectively, and suppose for simplicity that the supremum in Lemma~\ref{key0} is attained when $\Gamma=\{ (1,2),(1,3)\}$, say,  which we have identified with $\{2,3\}$. Then
$$
\|M_\Omega\|_{2\to 2} \, \le \ C\,\Big\| \sup_{{\bf i} \in \mathbb{Z}^{3}} M_{\Omega_{\bf
i}} K_{2, i_{2}}K_{3, i_{3}}\Big\|_{2\to
2}.
$$
Using the inclusion
$\ell^2(\mathbb{Z}^{2})\hookrightarrow
\ell^\infty(\mathbb{Z}^{2})$ and interchanging the order of the
sum and the integral,
\begin{align*}
\Big\| \sup_{{\bf i} \in \mathbb{Z}^{3}} M_{\Omega_{\bf i}} K_{2, i_{2}}K_{3, i_{3}}f\Big\|_2  &\ \le  \Big( \sum_{i_2,i_3 \in \mathbb{Z}} \big\| \sup_{i_1 \in \mathbb{Z}} M_{\Omega_{\bf i}}K_{2, i_{2}}K_{3, i_{3}}f\big\|_2^2 \Big)^{\frac12} \\ &\ \le\  \sup_{i_2,i_3 \in \mathbb{Z}}  \|\sup_{i_1\in \mathbb{Z}} M_{\Omega_{\bf i}}\|_{2 \to 2} \Big( \sum_{i_2,i_3 \in \mathbb{Z}} \big\| K_{2, i_{2}}K_{3, i_{3}}f \big\|_2^2 \Big)^{\frac12}\nonumber\\
&\ =\  \sup_{i_2,i_3 \in \mathbb{Z}}  \|M_{\Omega_{2,i_2}\cap\Omega_{3,i_3}}\|_{2 \to 2} \Big( \sum_{i_2,i_3 \in \mathbb{Z}} \big\| K_{2, i_{2}}K_{3, i_{3}}f \big\|_2^2 \Big)^{\frac12}\nonumber\\
&\ \le \  C\sup_{i_3 \in \mathbb{Z}}  \|M_{\Omega_{3,i_3}}\|_{2 \to 2} \Big( \sum_{i_3 \in \mathbb{Z}} \big\|K_{3, i_{3}}f \big\|_2^2 \Big)^{\frac12}\nonumber\\
&\ \le \  C\sup_{i_3 \in \mathbb{Z}}  \|M_{\Omega_{3,i_3}}\|_{2 \to 2}\|f \|_2,\nonumber
\end{align*}
and so we are done. In the final two inequalities we used nothing more that the finite overlapping of the two--dimensional conic frequency supports.
\vspace{1em}

More generally, we consider
$\mathbb{Z}^{\Sigma}=\mathbb{Z}^{\Gamma} \times \mathbb{Z}^{\Sigma
\setminus \Gamma}$, and given ${\bf i} = (i_\sigma)_{\sigma \in
\Sigma}$, we write ${\bf i}={\bf j}\times{\bf k}$ where  ${\bf
j}=(i_\sigma)_{\sigma \in \Gamma}$ and ${\bf k}=(i_\sigma)_{\sigma
\in \Sigma\backslash \Gamma}.$ Using the inclusion
$\ell^p(\mathbb{Z}^{\Gamma})\hookrightarrow
\ell^\infty(\mathbb{Z}^{\Gamma})$ and interchanging the order of the
sum and the integral,
\begin{align}\label{just}
\Big\| \sup_{{\bf i} \in \mathbb{Z}^{\Sigma}} M_{\Omega_{\bf i}}f_{\bf j} \Big\|_p  &\ \le  \Big( \sum_{{\bf j} \in \mathbb{Z}^{\Gamma}} \big\| \sup_{{\bf k} \in \mathbb{Z}^{\Sigma\backslash\Gamma}} M_{\Omega_{\bf i}}f_{\bf j} \big\|_p^p \Big)^{\frac1p} \\ &\ \le\  \sup_{{\bf j} \in \mathbb{Z}^{\Gamma}}  \|\sup_{{\bf k} \in \mathbb{Z}^{\Sigma\backslash\Gamma}} M_{\Omega_{\bf i}}\|_{p \to p} \Big( \sum_{{\bf j} \in \mathbb{Z}^{\Gamma}} \big\| f_{\bf j} \big\|_p^p \Big)^{\frac1p}\nonumber\\
&\ \le \  \sup_{\sigma \in \Sigma} \sup_{i\in \mathbb{Z}}
\|M_{\Omega_{\sigma\!,i}} \|_{p \to p} \Big\|\Big(\sum_{{\bf j} \in
\mathbb{Z}^{\Gamma}} |f_{\bf j}|^p\Big)^{\frac 1p}
\Big\|_p\nonumber.
\end{align}

\noindent\textbf{Case $p \ge 2$.} Using the inclusion $\ell^2(\mathbb{Z}^{\Gamma})\hookrightarrow
\ell^p(\mathbb{Z}^{\Gamma})$, from \eqref{just} we obtain
$$
\Big\| \sup_{{\bf i} \in \mathbb{Z}^{\Sigma}} M_{\Omega_{\bf i}}f_{\bf j} \Big\|_p  \ \le \  \sup_{\sigma \in \Sigma} \sup_{i\in \mathbb{Z}}
\|M_{\Omega_{\sigma\!,i}} \|_{p \to p} \Big\|\Big(\sum_{{\bf j} \in
\mathbb{Z}^{\Gamma}} |f_{\bf j}|^2\Big)^{\frac 12}
\Big\|_p.
$$
Taking $f_{\bf j}=\big[\prod_{\sigma \in \Gamma}
K_{\sigma\!,i_{\!\sigma}}\big]f$, where ${\bf j} =
(i_\sigma)_{\sigma \in \Gamma}$, and applying Lemmas~\ref{key0} and
\ref{square}, we obtain
the desired estimate.

\vspace{1em}

\noindent \textbf{Case $1 < p < 2$.} This is based on an argument of
M. Christ used in \cite{C,A} which refined the argument of
Nagel--Stein--Wainger~\cite{NSW}. We suppose initially that $\Omega$
is finite, so that by the triangle inequality and the
Hardy--Littlewood maximal theorem, $M_\Omega$ is bounded.
Then by interpolating between
$$
\big\| \sup_{{\bf i} \in \mathbb{Z}^{\Sigma}} M_{\Omega_{\bf i}}
f_{{\bf j}} \big\|_p \le \|M_{\Omega}\|_{p\to p}\big\|\sup_{{\bf j}
\in \mathbb{Z}^{\Gamma}} |f_{\bf j}| \,\big\|_p
$$
and \eqref{just},
we see that $\big\| \sup_{{\bf i} \in \mathbb{Z}^{\Sigma}} M_{\Omega_{\bf i}}
f_{\bf j}  \big\|_p$ is bounded above by
\begin{equation*} \label{Int}
 \|M_{\Omega}\|_{p\to p}^{1-\frac{p}{2}} \Big(
\sup_{\sigma \in \Sigma} \sup_{i \in\mathbb{Z}} \|
M_{\Omega_{\sigma\!,i}} \|_{p \to p} \Big)^{\frac{p}{2}} \Big\|\Big(
\sum_{{\bf j} \in \mathbb{Z}^{\Gamma}} |f_{\bf j}|^2 \Big)^\frac12
\Big\|_p.
\end{equation*}
Taking $f_{\bf j} = \big[\prod_{\sigma \in \Gamma}
K_{\sigma\!,i_{\!\sigma}}\big] f$, where~${\bf j} =
(i_\sigma)_{\sigma \in \Gamma}$, and applying Lemmas~\ref{key0} and
\ref{square} as before, we see that
$$
\|M_\Omega\|_{p\to p}\le C\,\|M_{\Omega}\|_{p\to p}^{1-\frac{p}{2}} \Big( \sup_{\sigma
\in \Sigma} \sup_{i \in\mathbb{Z}} \|M_{\Omega_{\sigma\!,i}}\|_{p\to p}
\Big)^{\frac{p}{2}}.
$$
Rearranging, we obtain the desired estimate with $C$ independent of $\Omega$, so we can drop the restriction that $\Omega$ is finite. This completes the proof. \fin

\vspace{1em}

In both \cite{NSW} and \cite{C},  a single conic Fourier multiplier
was introduced for each direction. This multiplier had to cover (the
bulk of) the hyperplane perpendicular to the direction, and so was
necessarily multidimensional in nature.  Rather restrictive conditions on the directions were then required to ensure finite overlapping of the supports of the
multipliers, yielding a bound via orthogonality as above.
In order to achieve greater flexibility, we introduced a number of essentially two--dimensional
multipliers instead. This is only possible via a covering rather than a partition, however after adding and subtracting a number of products of
these multipliers we obtain a {\it signed partition of unity}. This came at essentially no cost and in fact
simplifies matters because the orthogonality in two dimensions, summing over one index at a time, is
trivial to check. On the other hand, our multipliers are naturally associated to partitions of the directions allowing us to introduce a multiplier for
each segment instead of one for each direction.

\section{Proof of Theorem~C}\label{BB}

First we  prove the inclusion $\pmax(n)\subset \lacshadows(n)$
 which is restated in the following lemma. We appeal to a
quantitative version of Bateman's theorem~\cite{B}, allowing us to
treat the shadows simultaneously and thus uniformly. We also
use that the cross product of a two--dimensional Kakeya set with a
cube is a Kakeya set.

\begin{lemma}\label{bat} Let $n\ge 2$ and $1<p<\infty$, and suppose that $M_{\Omega}$ is bounded from $L^p(\mathbb{R}^n)$ to
$L^p(\mathbb{R}^n)$. Then $\Omega\in \lacshadows(n)$.
\end{lemma}

\dem As $M_\Omega$ is bounded if and only if $M_{\overline{\Omega}}$ is bounded, we can suppose that $\Omega$ is closed.
We appeal to Bateman's terminology \cite{B}. In particular we
will consider the binary tree $\mathcal{T}_{\Pi}$ associated to the
shadow of $\Omega$ on $\Pi$, for any two--dimensional subspace $\Pi$, and their splitting numbers
$\mathrm{split}(\mathcal{T}_{\Pi})$. We say that $\Omega$ {\it
admits Kakeya shadows} if there exists a constant $C$ such that
for any $N\ge 1$ there exists a two--dimensional subspace
${\Pi(N)}$ and a finite collection of rectangles
$\mathcal{R}_{\Pi(N)}$ contained in $\Pi(N)$, with longest side
pointing in a direction of the shadow of $\Omega$ on $\Pi(N)$, that
satisfy
\begin{equation}\label{theeq1}
\Big|\bigcup_{R\in \mathcal{R}_{\Pi(N)}}  R \ \Big| \le \frac{C}{N}
\, \Big| \bigcup_{R\in \mathcal{R}_{\Pi(N)}}  3R \ \Big|.
\end{equation}
Here, $3R$ has the same center and width as $R$, but three times the
length.

We prove the contrapositive. If $\Omega\not\in \lacshadows(n)$, then by Theorem 3 (combined with Remark 2) in \cite{B},
 for any $N\ge 1$, there is a shadow
of $\Omega$ on ${\Pi(N)}$ for which
$\mathrm{split}(\mathcal{T}_{\Pi(N)})\ge 2^N$. Bateman proved (see pages 61--62
 and Claim 7 of \cite{B}) that $\mathrm{split}(\mathcal{T}_{\Pi(N)}) \ge
2^N$ implies the existence of a finite family $\mathcal{R}_{\Pi(N)}$
of rectangles $R$ satisfying \eqref{theeq1}. Now for each
$N\in\mathbb{N}$, we pick an orthonormal basis $(e_1,\ldots,e_n)$ so
that $\text{span}(e_1,e_2)=\Pi(N)$. For each rectangle $R$ in the
subcollection $\mathcal{R}_{\Pi(N)}$, we set $$
\beta\equiv\beta(R)=\text{diam}(R)(\omega^2_1+\omega^2_2)^{-1/2},$$
where $\omega$ is a direction of $\Omega$ whose shadow points
in the direction of $R$, and let $\alpha\equiv\alpha(N)$ to be ten
times the maximum~$\beta(R)$ with~$R\in \mathcal{R}_{\Pi(N)}$.
Taking
$$E_N=\bigcup_{R\in \mathcal{R}_{\Pi(N)}}  R\times [0,\alpha]^{n-2},$$
defined with respect to the basis $(e_1,\ldots,e_n)$, we then have
$$M_{\Omega}[\chi_{E_N}](x)\ge 1/8\qquad  \text{for all}
\quad x\in \bigcup_{R\in \mathcal{R}_{\Pi(N)}}
 3R\times [3\beta,\alpha-3\beta]^{n-2}.$$
Using \eqref{theeq1}, we see that for all $N\ge 1$,
\begin{align*}
\big\|M_{\Omega}[\chi_{E_N}]\big\|_p\ge
cN^{\frac1p}\|\chi_{E_N}\|_p,
\end{align*}
so that $M_{\Omega}$ is not bounded from $L^p(\R^n)$ to $L^p(\R^n)$
when $p$ is finite.
 \fin

Using Theorem A in order to bound the maximal operators associated to the sets of  $\lac(n)$ and $\plac(n,\varepsilon)$, at this stage we have proven that
$$ \lac(n)\subset \plac(n)\subset\pmax(n)\subset \lacshadows(n).$$
It therefore remains to prove that $\pmax(n)\subset \plac(n)$. It is tempting to suppose that the job is already done -- that Theorem~A can be applied
repeatedly in order to reduce a shadow to a single direction,
thus reducing the dimension of the problem. That is to say $\lacshadows(n)\subset \lac(n)$ yielding a full chain of equivalences.
However the lacunary
orders of the shadows are unstable in the sense that shadows on
two--dimensional subspaces which are close can have dramatically
different lacunary orders, and so it is
 not clear that it helps to apply Theorem~A and then change the basis in
 order to apply it again. One may be faced each time with
 lacunary orders that are as bad as before, and the process may
 never end. This would not be a problem if a slightly more flexible version of Theorem~A were true, however the obvious candidates for such a theorem are false (see the following section).

 We get round the  problem by keeping track of the accumulation order, which is necessarily finite; see below. Although this purely topological quantity contains much less information than the lacunary orders, it has the advantage of being basis invariant. The key idea is to apply Theorem~A in such a way that if you are forced to change basis too often, losing track of the gain in the lacunary order of the shadows, at least you can be sure of a gain with respect to the  accumulation order. Indeed, by compactness we know that there are only a finite number of points of maximal accumulation order, and so by a finite splitting we can suppose that such a point is unique.
 Considering the three--dimensional case,  by positioning the basis vector $e_3$ at the unique point of maximal accumulation,  two thirds of the segments of a dissection no longer contain $e_3$ (see Figure 2) and in those cases we immediately obtain a reduction in accumulation order. Once the accumulation order has been reduced to zero,  we are left with only a finite number of directions by compactness, and so the associated operator is bounded. On the other hand, if we never obtain a reduction in  the accumulation order, then at each application of Theorem~A we have partitioned using the shadow orthogonal to $e_3$ (see Figure 3). In this case we simply rotate the directions around the $e_3$--axis at each application of the theorem -- $e_3$ is left unchanged -- and so we can be sure that we are reducing the lacunary order of the shadow orthogonal to $e_3$. Eventually we reduce to the case where the shadow is a single point and so we reduce to the two--dimensional problem.
 
  The reason we are able to carry through the details of this argument for $p$--lacunarity and not for lacunarity, is that we need only concern ourselves with a single dominating segment, avoiding problems regarding uniformity. More specifically, we perform finite splittings in order to reduce to the case where we have a unique point of maximal accumulation order. To perform this with each of the segments of an infinite partition, we would need a uniform bound on the the number of points of maximal accumulation order in each segment, and we do not know how to prove this. Indeed, we cannot be sure if this is a merely technical problem or if it could be reflected in the geometry of the directions.
 
 In order to pull--up this three--dimensional argument into higher dimensions we will need some more auxiliary definitions. Given an $m$--dimensional subspace $\Pi\subset \R^n$ and a set of directions $\Omega$, we
define the {\it $m$--shadow~$\Xi$ of $\Omega$ on $\Pi$} by
$$ \Xi=\Big\{ \frac{\proj_{\!\Pi}(\omega)}{|\proj_{\!\Pi}(\omega)|} \ : \ \omega \in \Omega \setminus \Pi^\perp \Big\} \, \subset \, \Pi\cap\mathbb{S}^{n-1},$$
where $\proj_{\!\Pi}$ denotes the orthogonal projection onto~$\Pi$ (see Figure 3  for an illustration with $n=3$ and $m=2$). Note that a $2$--shadow is the same thing as a shadow. 
For an $m$--shadow~$\Xi$, we consider
$A_k=\text{Ac}(A_{k-1})$, where $A_1=\text{Ac}(\Xi)$, the
accumulation points of~$\Xi$. We say that~$\Xi$ has {\it
accumulation order $K$} if $A_{K}$ is a finite set.

\begin{lemma}\label{acum} Let $2\le m\le n$ and $1<p<\infty$, and suppose that $M_{\Omega}$ is bounded from $L^p(\mathbb{R}^n)$ to
$L^p(\mathbb{R}^n)$. Then the $m$--shadows of $\Omega$ have uniformly bounded 
accumulation order.
\end{lemma}

\dem As $M_{\Omega}$ is bounded,   the 2--shadows are lacunary with uniform bounds on the lacunary orders by Lemma~\ref{bat}, so that in particular
the $2$--shadows of~$\Omega$ have uniformly bounded accumulation
order. Thus, it will suffice to prove that if the accumulation order
of  an $m$--shadow $\Xi$ of~$\Omega$ on~$\Pi$ is greater than $K$,
then there exists a $2$--shadow of $\Xi$, and hence also of
$\Omega$, whose accumulation order is greater than $K$. We take
$\xi\in A_{K+1}$ and consider a sequence $\{\xi_j\}_{j\ge1}$ in
$A_{K}$ which accumulates at $\xi$. Then for all but (at most) one
$(m-1)$--dimensional subspace of $\Pi$, the $(m-1)$--shadows of
$\{\xi_j\}_{j\ge1}$ on the $(m-1)$--dimensional subspaces accumulate
at the $(m-1)$--shadows of~$\xi$. Then we consider sequences in $A_{K-1}$
that accumulate at~$\xi_j$. Again for all but one
$(m-1)$--dimensional subspace of~$\Pi$, the $(m-1)$--shadows on the
$(m-1)$--dimensional subspaces accumulate at the $(m-1)$--shadows of $\xi_j$.
Continuing the process, we see that for all but a countable number
of $(m-1)$--dimensional subspaces of $\Pi$, the $(m-1)$--shadow of $\xi$ is
of accumulation order $\ge K+1$. We take one such shadow and repeat
the process. This yields an $(m-2)$--dimensional shadow of the
$(m-1)$--dimensional shadow of $\xi$, which is an
$(m-2)$--dimensional shadow of $\xi$, that is of accumulation order
$\ge K+1$. Repeating the process, we obtain the desired result.
 \fin

 As before we fix an auxiliary $\varepsilon>0$ and  say that an $m$--shadow~$\,\Xi\,$ of $\Omega$, with $\text{span}(\Xi)=\R^d$, is
 \begin{itemize}
 \item  {\it  $(n,p)$--lacunary of order~0}  if it consists of a single point
\item  
 {\it  $(n,p)$--lacunary of order $
L$}  if there are members $\{\Xi_{\sigma\!,i_{\!\sigma}}\}_{\sigma\in \Sigma(d)}$ of a dissection of $\Xi$
 which are
$(n,p)$--lacunary of order $\le L-1$ and for which the sets $\{\Omega_{\sigma\!,i_{\!\sigma}}\}_{\sigma\in \Sigma(d)}$ that shade them
 are dominating;
$$
\|M_{\Omega_{\sigma\!,i}}\|_{p\to p}\le
\|M_{\Omega_{\sigma\!,i_{\!\sigma}}}\|_{p\to p}+\varepsilon\quad \text{for
all}\quad i\in \mathbb{Z}^*.
$$
\end{itemize}
Note that in this context the dominating sets (which from now on we refer to as dominating segments) need only dominate the rest of their partition; not the whole dissection.
We say that an $m$--shadow $\Xi$ of $\Omega$ is {\it $(n,p)$--lacunary} if it is a finite union of sets which are  $(n,p)$--lacunary of finite order. 

The final inclusion, $\pmax(n)\subset \plac(n)$, is a consequence of the fact that \begin{equation}\label{pl}\pmax(n)\subset \plac(n,\varepsilon)\quad \text{ for all }\quad \varepsilon>0.\end{equation}
To see this,  we consider $\Omega\in \pmax(n)$. By Lemma~\ref{bat} we have that any subset of $\Omega$ is contained in $\lacshadows(n)$, so the 2--shadows of any subset of $\Omega$ are
$(n,p)$--lacunary (the existence of dominating segments is assured as the operator is bounded). Thus \eqref{pl} can be obtained by  $n-2$ applications of the following lemma, observing that if the $n$--shadow of a set is $(n,p)$--lacunary then the set is  $p$--lacunary. 

   \begin{figure}
\centering
 \begin{tikzpicture} 
\def\Ra{3} 
\def\angEl{115} 
\filldraw[ball color=white] (0,0) circle (\Ra);
\foreach \t in { 26.56505117, 45.,116.5650512,135.} { \DrawLongitudeCircle[\Ra]{\t} }
\pgfmathsetmacro\H{\Ra*cos(\angEl)}
\coordinate[mark coordinate] (N) at (0,\H);
\node[above=-15pt] at (N) {$e_{1}$};
\node[above=-9pt] at (N) {$\downarrow$};
\DrawLattwoCircle[\Ra]{0}
\coordinate[mark coordinate] (S) at (0,\H+3.95);
\node[above=12pt] at (S) {$e_{d}$};
\node[above=-0.5pt] at (S) {$\uparrow$};
\coordinate[mark coordinate] (M) at (0,\H);
\node[above=33pt, right=90 pt] at (M) {$\Pi$};
\node[above=-77pt] at (M) {{\bf Figure 2.} If $\sigma\in \Sigma(d)\backslash\Sigma(d-1)$ we separate from accumulation point.};
\coordinate[mark coordinate] (a) at (1.4,\H-1.2);
\coordinate[mark coordinate] (a) at (1.2,\H-1.2);
\coordinate[mark coordinate] (a) at (1,\H-1);
\coordinate[mark coordinate] (a) at (.8,\H-1);
\coordinate[mark coordinate] (a) at (.7,\H-0.8);

\coordinate[mark coordinate] (a) at (1.55,\H+1.3);
\coordinate[mark coordinate] (a) at (.75,\H+0.55);
\coordinate[mark coordinate] (a) at (.9,\H+0.76);
\coordinate[mark coordinate] (a) at (1.25,\H+0.85);
\coordinate[mark coordinate] (a) at (1.8,\H+1.85);
\coordinate[mark coordinate] (a) at (1.9,\H+1.76);
\coordinate[mark coordinate] (a) at (2.25,\H+1.85);
\coordinate[mark coordinate] (a) at (2.3,\H+2.85);
\coordinate[mark coordinate] (a) at (2.2,\H+2.76);
\coordinate[mark coordinate] (a) at (2.4,\H+2.65);

\coordinate[mark coordinate] (a) at (-2,\H+3.15);
\coordinate[mark coordinate] (a) at (-1.7,\H+3.45);
\coordinate[mark coordinate] (a) at (-1.7,\H+3);
\coordinate[mark coordinate] (a) at (-1.4,\H+3);
\coordinate[mark coordinate] (a) at (-1.4,\H+2.5);
\coordinate[mark coordinate] (a) at (-1.7,\H+2.5);
\coordinate[mark coordinate] (a) at (-1.3,\H+2);
\coordinate[mark coordinate] (a) at (-1,\H+1.7);
\coordinate[mark coordinate] (a) at (-.9,\H+1.3);
\coordinate[mark coordinate] (a) at (-.6,\H+.9);

\coordinate[mark coordinate] (a) at (-.8,\H+-.5);
\coordinate[mark coordinate] (a) at (-1.1,\H+-.7);
\coordinate[mark coordinate] (a) at (-1.3,\H+-.8);
\coordinate[mark coordinate] (a) at (-1.6,\H+-.8);
\coordinate[mark coordinate] (a) at (-1.9,\H+-.8);
\coordinate[mark coordinate] (a) at (-2.2,\H+-.6);

\end{tikzpicture}
\end{figure}

\begin{lemma}\label{final} Let $2\le m\le n-1$ and $1<p<\infty$, and suppose that $M_{\Omega}$ is bounded from $L^p(\R^n)$ to $L^p(\R^n)$. Then, if the $m$--shadows of any subset of $\Omega$ are $(n,p)$--lacunary, then the $(m+1)$--shadows of any subset of
$\Omega$ are $(n,p)$--lacunary.
\end{lemma}

\dem 
Consider the $(m+1)$--shadow of an arbitrary subset  $\Omega_0\subset\Omega\in \pmax(n)$. By compactness there are a finite number of points of maximal accumulation order, and after a finite splitting we can suppose that each piece has a unique point of maximal accumulation, and that the $m$--shadow orthogonal to this is $(n,p)$--lacunary of finite order.  Here we use that the $m$--shadow 
of an  $(m+1)$--shadow of $\Omega_0$ is the same as the $m$--shadow of $\Omega_0$, and so it is $(n,p)$--lacunary by hypothesis.
It will suffice to prove that each of these pieces of the  $(m+1)$--shadow, we consider an arbitrary piece $\Xi$, is  $(n,p)$--lacunary of finite order.

  If the accumulation order of $\Xi$ were zero, then it would consist of a single point, and so it would be $(n,p)$--lacunary of order zero.  On the other hand, if the $m$--shadow, orthogonal to the point of maximal accumulation order, were of $(n,p)$--lacunary order zero, then it would also consist of a single point. In this case,  $\Xi$ would be contained in a two--dimensional
subspace, and so $\Xi$ would form part of a $2$--shadow of $\Omega$. From here we could deduce that $\Xi$ would again be $(n,p)$--lacunary of finite order by
Lemma~\ref{bat}.

 \begin{figure}
\centering
 \begin{tikzpicture} 
\def\Ra{3} 
\def\angEl{25} 
\filldraw[ball color=white] (0,0) circle (\Ra);
\foreach \t in  { 26.56505117, 45.,116.5650512,135.}  { \DrawLongitudeCircle[\Ra]{\t} }
\foreach \t in  {30.5,34,41,118.5,122.5,130,135.}  { \DrawLongitudeCircletwo[\Ra]{\t} }
\DrawLatitudeCircle[\Ra]{0}
\pgfmathsetmacro\H{\Ra*cos(\angEl)}
\coordinate[mark coordinate] (N) at (0,\H);
\node[above=11pt] at (N) {$e_{d}$};
\node[above=-2pt] at (N) {$\uparrow$};
\coordinate[mark coordinate] (M) at (0,\H);
\node[above=-80pt, right=90 pt] at (M) {$\Pi$};
\node[above=-190pt] at (M) {{\bf Figure 3.} If $\sigma\in \Sigma(d-1)$ we reduce the order of the shadow on $\Pi$.};
\coordinate[mark coordinate] (a) at (1.87,\H-3);
\coordinate[mark coordinate] (a) at (1.92,\H-4);
\coordinate[mark coordinate] (b) at (1.12,\H-1.5);
\coordinate[mark coordinate] (b) at (.97,\H-1.2);
\coordinate[mark coordinate] (b) at (1.07,\H-1);
\coordinate[mark coordinate] (b) at (.54,\H-0.5);
\coordinate[mark coordinate] (b) at (.44,\H-0.3);
\coordinate[mark coordinate] (b) at (.23,\H-0.15);
\coordinate[mark coordinate] (b) at (1.27,\H-5);
\coordinate[mark coordinate] (b) at (1.48,\H-4.8);
\coordinate[mark coordinate] (b) at (1.4,\H-3.3);
\coordinate[mark coordinate] (b) at (1.51,\H-2.8);
\coordinate[mark coordinate] (b) at (-.6,\H-.3);
\coordinate[mark coordinate] (a) at (-2.22,\H-3);
\coordinate[mark coordinate] (a) at (-2.12,\H-1.7);
\coordinate[mark coordinate] (a) at (-2.1,\H-1.87);
\coordinate[mark coordinate] (a) at (-1.68,\H-1.2);
\coordinate[mark coordinate] (a) at (-1.6,\H-1.37);
\coordinate[mark coordinate] (a) at (-1.5,\H-1);
\coordinate[mark coordinate] (a) at (-1.5,\H-.85);
\coordinate[mark coordinate] (a) at (-1,\H-.65);
\coordinate[mark coordinate] (a) at (-1,\H-.45);
\coordinate[mark coordinate] (a) at (-.4,\H-.15);
\coordinate[mark coordinate] (b) at (-2.25,\H-4);
\coordinate[mark coordinate] (b) at (-2.26,\H-4.5);
\coordinate[mark coordinate] (b) at (-0.7,\H-.3);
\coordinate[mark coordinate] (b) at (0.9,\H+.09);
\end{tikzpicture}
\end{figure}

Thus we know that $\Xi$ is $(n,p)$--lacunary of finite order if the accumulation order $K$ is zero or if the $(n,p)$--lacunary order $L$ of  the  $m$--shadow, orthogonal to the point of maximal accumulation order, is zero. 
Supposing that we also knew this to be true if the accumulation order is $\le K$ or if the orthogonal $m$--shadow is of $(n,p)$--lacunary order $\le L$, by induction it would suffice to prove that $\Xi$ is $(n,p)$--lacunary of finite order supposing that the accumulation order is $K+1$ and the $m$--shadow of $\Xi$, orthogonal to the point of maximal accumulation order, is of $(n,p)$--lacunary order $L+1$.

In order to prove this induction step, we carefully choose the basis vectors $e_1,\ldots,e_{d}$ in order to dissect $\Xi$, where $d\le m+1$ is the dimension of $\text{span}(\Xi)$. It would suffice to find a dissection for which all of the dominating segments are $(n,p)$--lacunary of finite order, as then we could conclude that $\Xi$ is $(n,p)$--lacunary of finite order. We take $e_{d}$ to be the point of maximal accumulation order and choose the remaining vectors in order to  dissect $\Xi$ (simultaneously dissecting the
$m$--shadow that lives in~$\Pi=\text{span}(e_1,\ldots,e_{d-1})$ and partially dissecting $\Omega$) with
$(e_1,\ldots,e_{d-1})$ and
 $\{\theta_{\sigma\!,i}\}_{i\in\mathbb{Z}}$ for each $\sigma\in \Sigma(d-1)$ chosen in order
 to reduce the $(n,p)$--lacunary order of the $m$--shadow that lives in $\Pi$
  (see Figure 3). We are free to choose any lacunary $\{\theta_{\sigma\!,i}\}_{i\in\mathbb{Z}}$ for $\sigma=(j,d)$
  with $1\le j\le d-1$. There are dominating segments in each partition as $M_\Omega$ is bounded. Those of the partitions with $\sigma=(j,d)$ are separated from $e_{d}$ (see Figure 2), so they have reduced  accumulation order $\le K$, and so they are $(n,p)$--lacunary of finite order by the induction hypothesis. On the other hand, the $m$--shadows that live in $\Pi$ of the dominating segments of the other partitions have reduced $(n,p)$--lacunary order $\le L$, and so these dominating segments are also $(n,p)$--lacunary of finite order by the induction hypothesis. Altogether, we have found a dissection of  $\Xi$ with dominating segments which are all $(n,p)$--lacunary of finite order.  Thus  $\Xi$ is $(n,p)$--lacunary of finite order, and the proof is complete by induction. \fin

\section{Our notion of lacunarity and the sharpness of Theorem A}\label{three3}

As in the
previous section, to construct unbounded directional maximal operators, the directions need only be badly spaced after
projecting onto a two--dimensional subspace. Thus, in contrast
with the two--dimensional case, it is not
 enough  to constrain the angles
 between the directions if they are to give rise to a bounded
 maximal operator in higher dimensions.
To see this, we enumerate  $\mathbb{Q}\cap[\frac{1}{2},\frac{2}{3}]=\{q_\ell\}_{\ell\ge 1}$ and consider
\begin{equation*}
\Omega=\Big\{\, \omega \in \mathbb{S}^{n-1}\cap\R_+^n\, :\,  \frac{\omega_2}{\omega_1}=q_\ell,\  \omega_j=2^{-j\ell},   1< j<n;\ \text{for some}\
   \ell\ge
1\,\Big\}.\end{equation*}
Then the angles between the directions form a lacunary sequence converging to zero with lacunary constant $1/2$; see Figure 4. Taking~$\theta_{\sigma\!,i}=2^{-i}$, the segments $\Omega_{\sigma\!,i}$
consist of at most one direction for all $i\in \mathbb{Z}^*$ and
$\sigma\in \Sigma\backslash \{(1,2)\}$.
In spite of this, $M_\Omega$
 is unbounded. Indeed, consider the set of
rectangles $\mathcal{R}$ in $\Pi=\text{span}(e_1,e_2)$ with longest
side parallel to the shadow on $\Pi$ of some~$\omega\in \Omega$. Then
the construction of Besicovitch (see for example \cite{F}) provides finite
subsets $\mathcal{R}_N\subset \mathcal{R}$, for all~$N\ge 1$,  that
satisfy \eqref{theeq1}. Considering $\chi_{E_N}$, defined as in the
proof of Lemma~\ref{bat}, we find $M_{\Omega}$ unbounded as
before.

If the angles between directions restricted to a great circle are lacunary, or if the angles between  directions restricted to the Nagel--Stein--Wainger curves are lacunary, then the associated maximal operators are bounded.
It is tempting to suppose that if the angles between directions restricted to any smooth curve (which does not spiral around  the sphere infinitely many times) are lacunary then the directions give rise to a bounded maximal operator (the authors thank Antonio C\'ordoba for asking this question). To see that this is not the case we consider the curve $\gamma:[0,1/4] \to \mathbb{S}^{n-1}$ defined to be the normalisation of $\widetilde{\gamma}(t)=(t,t/\log_2(1/t),t,\ldots,t,1)$. This is little more than a smooth perturbation of a great circle (see also \cite{MS} for a counterexample for the directional maximal operator defined on manifolds with degenerate curvature). We consider the directions $\Omega=\{\omega_\ell\}_{\ell\ge 1}$ where $\omega_\ell=\gamma(2^{-\ell})$; see Figure 5. As long as $\ell$ is taken sufficiently large we can safely ignore the normalisation. 
Then it is easy to see that the angles between the directions are lacunary with lacunary constant $1/2$:
\begin{align*}
\frac{|\widetilde{\gamma}(2^{-(\ell+1)})-e_n|}{|\widetilde{\gamma}(2^{-\ell})-e_n|}&\le \frac{|(2^{-(\ell+1)},2^{-(\ell+1)}/(\ell+1),\ldots,2^{-(\ell+1)},0)|}{|(2^{-\ell},2^{-\ell}/\ell,\ldots,2^{-\ell},0)|}\\
&\le \frac{1}{2}\frac{|(1,1/(\ell+1),1,\ldots,1,0)|}{|(1,1/\ell,1,\ldots,1,0)|}< \frac{1}{2}.
\end{align*}
In spite of this, $M_\Omega$
 is unbounded. Indeed, consider the set of
rectangles $\mathcal{R}$ in $\Pi=\text{span}(e_1,e_2)$ with longest
side parallel to the shadow on $\Pi$ of some~$\omega\in \Omega$. Then
the construction of Besicovitch provides finite
subsets $\mathcal{R}_N\subset \mathcal{R}$, for all~$N\ge 1$,  that
satisfy \eqref{theeq1}. 
 To see this it is enough to show that there are approximately uniformly spaced  angles between the shadows of the directions at all scales. We have that
 $$
 \frac{\omega_\ell\cdot e_2}{\omega_\ell\cdot e_1}= \frac{1}{\ell}\qquad\textrm{and}\qquad \frac{\omega_\ell\cdot e_2}{\omega_\ell\cdot e_1}-\frac{\omega_{\ell+1}\cdot e_2}{\omega_{\ell+1}\cdot e_1}=\frac{1}{\ell(\ell+1)}
 $$
 so that the shadow contains $\ell$ approximately equally spaced points between the shadow of $\omega_\ell$ and $e_1$ for all $\ell$ sufficiently large.
Considering $\chi_{E_N}$, defined as in the
proof of Lemma~\ref{bat}, we find $M_{\Omega}$ unbounded as
before.

\begin{figure}
\centering
\begin{tikzpicture} 
\def\Ra{3} 
\def\angEl{5} 
\filldraw[ball color=white] (0,0) circle (\Ra);
\foreach \t in  {116, 118, 120, 122, 124, 126, 128, 130, 132, 134}  { \DrawLongCircle[\Ra]{\t} }
\DrawLatitudeCircle[\Ra]{0}
\pgfmathsetmacro\H{\Ra*cos(\angEl)}
\coordinate[mark coordinate] (N) at (0,\H);
\node[above=10pt] at (N) {$e_{n}$};
\node[above=-1pt] at (N) {$\uparrow$};
\coordinate[mark coordinate] (a) at (2,\H-3);
\coordinate[mark coordinate] (b) at (1.17,\H-1.5);
\coordinate[mark coordinate] (b) at (.8,\H-.75);
\coordinate[mark coordinate] (b) at (.8,\H-.3);
\coordinate[mark coordinate] (b) at (.4,\H-.18);
\coordinate[mark coordinate] (b) at (.22,\H-.065);
\coordinate[mark coordinate] (b) at (.1,\H-.02);
\node[above=-85pt, right=90 pt] at (M) {$\Pi$};
\coordinate[mark coordinate] (M) at (0,\H);
\node[above=-200pt] at (M) {{\bf Figure 4.} A Kakeya shadow on $\Pi$ with $\frac{1}{2}n(n-1)-1$ lacunary shadows.};
\end{tikzpicture}

\end{figure}

 If Theorem A were more flexible, in the sense that the partitions were allowed to \lq accumulate' away from the hyperplanes orthogonal to the basis vectors, then we would obtain $\pmax(n)\equiv\lacshadows(n)$ as we would be able to bound the
operators associated to the sets of $\lacshadows(n)$. However, Theorem~A is remarkably
sharp in the sense that the supremum in $\sigma$ must be taken over
the whole of~$\Sigma$, and the partitions must accumulate at the hyperplanes perpendicular to the basis vectors. To see this, we let $e_2'$ and $e'_n$ be orthogonal unit vectors in $\text{span} (e_2,e_n)$, close to $e_2$ and $e_n$, with $e_n'$ in the first quadrant determined by $e_2$ and~$e_n$. We construct a set of directions, accumulating rapidly at $e_n'$,
for which the angles between the orthogonal projections onto
$\text{span}(e_1,e_{2}')$ are badly spaced. Indeed, we  take
$\Omega=\{\omega_\ell\}_{\ell\ge 1}$ so that $\omega_\ell\cdot e_2'=q_\ell\,
\omega_\ell\cdot e_1$. This does not yet completely determine~$\omega_\ell$.
Supposing that we have chosen $\omega_{\ell-1}$ we can choose the
direction~$\omega_\ell$ sufficiently close to $e_n'$ so that the
angle between $\omega_{\ell-1}$ and $e_{n}'$ is at least double that
between $\omega_{\ell}$ and~$e_n'$. We can also choose the directions so that
$$
\frac{\omega_{\ell-1}\cdot e'_n}{\omega_{\ell-1}\cdot e'_2} \le \frac12 \frac{\omega_{\ell}\cdot e'_n}{\omega_{\ell}\cdot e'_2},\qquad
\text{and}\qquad
\frac{\omega_{\ell-1}\cdot e_k}{\omega_{\ell-1}\cdot e_j} \le \frac12 \frac{\omega_{\ell}\cdot e_k}{\omega_{\ell}\cdot e_j}
$$
for all $(j,k)\in\Sigma(n)\backslash \{(2,n)\}.$
  Taking~$\theta_{\sigma\!,i}=2^{-i}$, the segments
  $\Omega_{\sigma\!,i}$, defined with respect to the orthonormal basis $(e_1,\ldots,e_n)$,
consist of at most one direction for all~$i\in \mathbb{Z}^*$
and~$\sigma\in \Sigma(n)\backslash \{(2,n)\}$. On the other hand, if we define the final segments by
$$ \Omega_{(2,n),i}=\Big\{\, \omega \in \Omega \ :\ 2^{-(i+1)}< \Big|\frac{\omega\cdot e'_n}{\omega\cdot e_2'}\Big|\le 2^{-i}\, \Big\},\quad i\in\mathbb{Z},$$
accumulating at $\{e'_2\}^\perp\cup \{e'_n\}^\perp,$
then they also consist of at most one direction for all ~$i\in \mathbb{Z}$.
In spite of this, $M_\Omega$
 is unbounded as before. Indeed, consider the set of
rectangles $\mathcal{R}$ in $\Pi=\text{span}(e_1,e_2')$ with longest
side parallel to the shadow on $\Pi$ of some~$\omega_\ell$.
Then there are finite subsets $\mathcal{R}_N\subset \mathcal{R}$,
for all~$N\ge 1$,  that satisfy \eqref{theeq1}. Considering
$\chi_{E_N}$, defined as in the proof of Lemma~\ref{bat}, but with
respect to the basis $(e_1,e_2',e_3,\ldots,e_{n-1},e'_n)$, we again find
$M_{\Omega}$ unbounded on $L^p(\R^n)$ for finite~$p$.

\begin{figure}
\centering
\begin{tikzpicture} 
\def\Ra{3} 
\def\angEl{5} 
\filldraw[ball color=white] (0,0) circle (\Ra);
\foreach \t in  {70, 68, 66, 64, 62, 60}  { \DrawLongtwoCircle[\Ra]{\t} }
\DrawLatitudeCircle[\Ra]{0}
\pgfmathsetmacro\H{\Ra*cos(\angEl)}
\coordinate[mark coordinate] (N) at (0,\H);
\node[above=10pt] at (N) {$e_{n}$};
\node[above=-1pt] at (N) {$\uparrow$};
\coordinate[mark coordinate] (b) at (-.84,\H-1.5);
\coordinate[mark coordinate] (b) at (-.68,\H-.75);
\coordinate[mark coordinate] (b) at (-.45,\H-.3);
\coordinate[mark coordinate] (b) at (-.28,\H-.13);
\node[above=-92pt, left=39pt] at (N) {$\leftarrow$};
\node[above=-92pt, left=50pt] at (N) {$e_1$};
\coordinate[mark coordinate] (b) at (-1.52,\H-3.2);
\coordinate[mark coordinate] (b) at (-.13,\H-.05);
\node[above=-85pt, right=90 pt] at (M) {$\Pi$};
\coordinate[mark coordinate] (M) at (0,\H);
\node[above=-195pt] at (M) {{\bf Figure 5.} A Kakeya shadow of well-spaced directions in a smooth curve.};
\end{tikzpicture}
\end{figure}

Finally we remark that \lq cross products' of lacunary sets, like the directions~\eqref{carbery}, do not give rise to bounded maximal operators in general. To see this we consider the largest set of the form
\begin{equation*}
\Omega=\Big\{ \omega \in \mathbb{S}^{n-1}\cap\R_+^n:  \frac{\omega_k}{\omega_j}=2^{-i},\  \frac{\omega_n}{\omega_1}=3^{-\ell},\   1< j<k;\ \text{for some}\
   i,\ell\in \mathbb{Z}\Big\},\end{equation*}
and $\Theta=\{2^{i}3^{-\ell}\}_{i,\ell\in \mathbb{Z}}$ which is the set of the tangents of the angles between the shadows of the directions on $\Pi=\text{span}(e_1,e_2)$. To see that this is dense in $\R_+$, which is presumably well-known, we note that 
$$
|2^{i}3^{-\ell}-1|< \epsilon\ \Leftrightarrow\ \big|\tfrac{i}{\ell}-\log_2 3\big|<\tfrac{\log_2(1+\epsilon)}{\ell},
$$
when $2^{i}3^{-\ell}>1$, so that by Dirichlet's approximation theorem, 1 is an accumulation point of $\Theta$. Then if $\Theta$ were not dense we could find an interval $(a,b)$, with $a,b$ in the closure of $\Theta$, which does not contain an element of~$\Theta$. However, noting that $\Theta$ is closed under multiplication, by taking a sequence of $\Theta$ which accumulates to $1$ from above and multiplying by elements of $\Theta$ sufficiently close to $a$, we come to a contradiction. Considering $\chi_{E_N}$, defined as in the
proof of Lemma~\ref{bat}, we find $M_{\Omega}$ unbounded as
before.

\section{The maximal directional Hilbert transform}

It is well known that there is a close relationship between the behaviour of the directional maximal operator and the maximal directional Hilbert transform $H_\Omega$, defined by
$f\mapsto \sup_{\omega\in \Omega} |H_\omega f|,
$
where
$$
H_\omega f(x)= p.v. \int_{\R} f(x-\omega t)\,\frac{dt}{t}.
$$
However, the constant in the following corollary must depend on the cardinality of $\Omega$ due to a result of Karagulyan  which showed that the maximal directional Hilbert transform in the plane is unbounded as soon as the number of directions is infinite \cite{Ka}. On the other hand, we do not recover the sharp estimates for $H_\Omega$ in terms of the power of the logarithm when $n=2$; see \cite{DdP},  and so it would be interesting to see if the following inequality could be improved in that regard. The estimate also holds for more general operators, where the kernel $1/t$  is replaced by the inverse Fourier transform of a H\"ormander--Mikhlin multiplier. 
\begin{corollary} Let $n\ge 2$ and $p>1$. Then
\begin{align*}
 \|H_{\Omega}f\|_{p\to p}\ \le\  C\log |\Omega|\, \sup_{\sigma \in \Sigma}\, \sup_{i\in\mathbb{Z}^*}\|M_{\Omega_{\sigma\!,i}}\|_{p\to p},
\end{align*}
where $C$ depends only on $n$, $p$ and the lacunary constants $\lambda_\sigma$ for $\sigma \in \Sigma$.
\end{corollary}

\begin{proof}  First we note that $\{ x\,:\, H_{\Omega}f(x)>\log |\Omega|\gamma\}$ is a subset of 
\begin{align*}
 \big\{x\,:\, H_{\Omega}f(x)>\log |\Omega|\gamma,\ M_\Omega f(x) \le \gamma\big\}\cup\big\{x\,:\,  M_{\Omega}f(x)>\gamma\big\},
\end{align*}
so we can use Theorem A to deal with the part of the integral coming from the second level set. Thus,  it will suffice to prove
\begin{equation}\label{pok}
\int_0^\infty \!\big|\big\{ x\, :\, H_{\Omega}f(x)>\log |\Omega|\gamma,\, M_\Omega f(x)  \le  \gamma\big\}\big|p\gamma^{p-1} d\gamma \le  C\|f\|_p^p.
\end{equation}
 To see this we first note that
 \begin{align*}
 & \ \big|\big\{ x\, :\, H_{\Omega}f(x)>\log |\Omega|\gamma,\ M_\Omega f(x) \ \le \ \gamma\big\}\big|\\\le\ &\sum_{\omega\in \Omega}\big|\big\{ x\, :\, H_{\omega}f(x)>\log |\Omega|\gamma,\ M_\omega f(x) \ \le \ \gamma\big\}\big|.
 \end{align*}
 Then we use a reformulation of a  one--dimensional inequality due to Hunt, 
 $$
 \big|\big\{ x\, :\, H_{\omega}f(x)>N\gamma,\ M_\omega f(x) \ \le \ \gamma\big\}\big|\le  e^{-N} \big|\big\{ x\, :\, H^\star_{\omega}f(x)>\gamma \big\}\big|
 $$
 (see \cite[Proposition 2.2]{DdP} for more details),
 where
 $$
 H^\star_\omega f(x)= \sup_{\epsilon>0}\Big|\int_{|t|>\epsilon} f(x-\omega t)\,\frac{dt}{t}\Big|.
 $$
 Altogether we see that the left-hand side of \eqref{pok} is bounded by a constant multiple of
 \begin{align*}
\ &\ \sum_{\omega\in \Omega} \frac{1}{|\Omega|}\int_0^\infty\big|\big\{ x\, :\, H^\star_{\omega}f(x)>\gamma \big\}\big|p\gamma^{p-1} d\gamma\\
\le &\  \sum_{\omega\in \Omega} \frac{1}{|\Omega|}{\|H^\star_{\omega}f\|_p^p}\,\le\, C\sum_{\omega\in \Omega}\frac{1}{|\Omega|}{\|f\|_p^p}\,\le\, C{\|f\|_p^p},
 \end{align*}
and so we are done.
\end{proof}

\noindent The authors thank Francesco Di Plinio for pointing out the final corollary and the referees for helpful comments.

\end{document}